\documentclass{amsart}

\input{epsf.tex}
\usepackage{amsfonts,amssymb,verbatim,amsmath,amsthm,latexsym,textcomp,amscd}
\usepackage{latexsym,amsfonts,amssymb,epsfig,verbatim}
\usepackage{amsmath,amsthm,amssymb,latexsym,graphics,textcomp}
\usepackage{graphicx}
\usepackage{color}
\usepackage{url}
\input{psfig.sty}

\begin{document}

\newtheorem{theorem}{Theorem}[section]
\newtheorem{prop}[theorem]{Proposition}
\newtheorem{lemma}[theorem]{Lemma}
\newtheorem{cor}[theorem]{Corollary}
\newtheorem{defn}[theorem]{Definition}
\newtheorem{conj}[theorem]{Conjecture}
\newtheorem{claim}[theorem]{Claim}

\newcommand{\boundary}{\partial}
\newcommand{\C}{{\mathbb C}}
\newcommand{\integers}{{\mathbb Z}}
\newcommand{\natls}{{\mathbb N}}
\newcommand{\ratls}{{\mathbb Q}}
\newcommand{\reals}{{\mathbb R}}
\newcommand{\proj}{{\mathbb P}}
\newcommand{\lhp}{{\mathbb L}}
\newcommand{\tube}{{\mathbb T}}
\newcommand{\cusp}{{\mathbb P}}
\newcommand\AAA{{\mathcal A}}
\newcommand\BB{{\mathcal B}}
\newcommand\CC{{\mathcal C}}
\newcommand\DD{{\mathcal D}}
\newcommand\EE{{\mathcal E}}
\newcommand\FF{{\mathcal F}}
\newcommand\GG{{\mathcal G}}
\newcommand\HH{{\mathcal H}}
\newcommand\II{{\mathcal I}}
\newcommand\JJ{{\mathcal J}}
\newcommand\KK{{\mathcal K}}
\newcommand\LL{{\mathcal L}}
\newcommand\MM{{\mathcal M}}
\newcommand\NN{{\mathcal N}}
\newcommand\OO{{\mathcal O}}
\newcommand\PP{{\mathcal P}}
\newcommand\QQ{{\mathcal Q}}
\newcommand\RR{{\mathcal R}}
\newcommand\SSS{{\mathcal S}}
\newcommand\TT{{\mathcal T}}
\newcommand\UU{{\mathcal U}}
\newcommand\VV{{\mathcal V}}
\newcommand\WW{{\mathcal W}}
\newcommand\XX{{\mathcal X}}
\newcommand\YY{{\mathcal Y}}
\newcommand\ZZ{{\mathcal Z}}
\newcommand\CH{{\CC\HH}}
\newcommand\TC{{\TT\CC}}
\newcommand\EXH{{ \EE (X, \HH )}}
\newcommand\GXH{{ \GG (X, \HH )}}
\newcommand\GYH{{ \GG (Y, \HH )}}
\newcommand\PEX{{\PP\EE  (X, \HH , \GG , \LL )}}
\newcommand\MF{{\MM\FF}}
\newcommand\PMF{{\PP\kern-2pt\MM\FF}}
\newcommand\ML{{\MM\LL}}
\newcommand\PML{{\PP\kern-2pt\MM\LL}}
\newcommand\GL{{\GG\LL}}
\newcommand\Pol{{\mathcal P}}
\newcommand\half{{\textstyle{\frac12}}}
\newcommand\Half{{\frac12}}
\newcommand\Mod{\operatorname{Mod}}
\newcommand\Area{\operatorname{Area}}
\newcommand\ep{\epsilon}
\newcommand\hhat{\widehat}
\newcommand\Proj{{\mathbf P}}
\newcommand\U{{\mathbf U}}
 \newcommand\Hyp{{\mathbf H}}
\newcommand\D{{\mathbf D}}
\newcommand\Z{{\mathbb Z}}
\newcommand\R{{\mathbb R}}
\newcommand\Q{{\mathbb Q}}
\newcommand\E{{\mathbb E}}
\newcommand\til{\widetilde}
\newcommand\length{\operatorname{length}}
\newcommand\tr{\operatorname{tr}}
\newcommand\gesim{\succ}
\newcommand\lesim{\prec}
\newcommand\simle{\lesim}
\newcommand\simge{\gesim}
\newcommand{\simmult}{\asymp}
\newcommand{\simadd}{\mathrel{\overset{\text{\tiny $+$}}{\sim}}}
\newcommand{\ssm}{\setminus}
\newcommand{\diam}{\operatorname{diam}}
\newcommand{\pair}[1]{\langle #1\rangle}
\newcommand{\T}{{\mathbf T}}
\newcommand{\inj}{\operatorname{inj}}
\newcommand{\pleat}{\operatorname{\mathbf{pleat}}}
\newcommand{\short}{\operatorname{\mathbf{short}}}
\newcommand{\vertices}{\operatorname{vert}}
\newcommand{\collar}{\operatorname{\mathbf{collar}}}
\newcommand{\bcollar}{\operatorname{\overline{\mathbf{collar}}}}
\newcommand{\I}{{\mathbf I}}
\newcommand{\tprec}{\prec_t}
\newcommand{\fprec}{\prec_f}
\newcommand{\bprec}{\prec_b}
\newcommand{\pprec}{\prec_p}
\newcommand{\ppreceq}{\preceq_p}
\newcommand{\sprec}{\prec_s}
\newcommand{\cpreceq}{\preceq_c}
\newcommand{\cprec}{\prec_c}
\newcommand{\topprec}{\prec_{\rm top}}
\newcommand{\Topprec}{\prec_{\rm TOP}}
\newcommand{\fsub}{\mathrel{\scriptstyle\searrow}}
\newcommand{\bsub}{\mathrel{\scriptstyle\swarrow}}
\newcommand{\fsubd}{\mathrel{{\scriptstyle\searrow}\kern-1ex^d\kern0.5ex}}
\newcommand{\bsubd}{\mathrel{{\scriptstyle\swarrow}\kern-1.6ex^d\kern0.8ex}}
\newcommand{\fsubeq}{\mathrel{\raise-.7ex\hbox{$\overset{\searrow}{=}$}}}
\newcommand{\bsubeq}{\mathrel{\raise-.7ex\hbox{$\overset{\swarrow}{=}$}}}
\newcommand{\tw}{\operatorname{tw}}
\newcommand{\base}{\operatorname{base}}
\newcommand{\trans}{\operatorname{trans}}
\newcommand{\rest}{|_}
\newcommand{\bbar}{\overline}
\newcommand{\UML}{\operatorname{\UU\MM\LL}}
\newcommand{\EL}{\mathcal{EL}}
\newcommand{\tsum}{\sideset{}{'}\sum}
\newcommand{\tsh}[1]{\left\{\kern-.9ex\left\{#1\right\}\kern-.9ex\right\}}
\newcommand{\Tsh}[2]{\tsh{#2}_{#1}}
\newcommand{\qeq}{\mathrel{\approx}}
\newcommand{\Qeq}[1]{\mathrel{\approx_{#1}}}
\newcommand{\qle}{\lesssim}
\newcommand{\Qle}[1]{\mathrel{\lesssim_{#1}}}
\newcommand{\simp}{\operatorname{simp}}
\newcommand{\vsucc}{\operatorname{succ}}
\newcommand{\vpred}{\operatorname{pred}}
\newcommand\fhalf[1]{\overrightarrow {#1}}
\newcommand\bhalf[1]{\overleftarrow {#1}}
\newcommand\sleft{_{\text{left}}}
\newcommand\sright{_{\text{right}}}
\newcommand\sbtop{_{\text{top}}}
\newcommand\sbot{_{\text{bot}}}
\newcommand\sll{_{\mathbf l}}
\newcommand\srr{_{\mathbf r}}
\newcommand\geod{\operatorname{\mathbf g}}
\newcommand\mtorus[1]{\boundary U(#1)}
\newcommand\A{\mathbf A}
\newcommand\Aleft[1]{\A\sleft(#1)}
\newcommand\Aright[1]{\A\sright(#1)}
\newcommand\Atop[1]{\A\sbtop(#1)}
\newcommand\Abot[1]{\A\sbot(#1)}
\newcommand\boundvert{{\boundary_{||}}}
\newcommand\storus[1]{U(#1)}
\newcommand\Momega{\omega_M}
\newcommand\nomega{\omega_\nu}
\newcommand\twist{\operatorname{tw}}
\newcommand\modl{M_\nu}
\newcommand\MT{{\mathbb T}}
\newcommand\Teich{{\mathcal T}}
\renewcommand{\Re}{\operatorname{Re}}
\renewcommand{\Im}{\operatorname{Im}}

\title{Flows, Fixed Points and Rigidity for Kleinian Groups}

\author{Kingshook Biswas }

\date{}

\maketitle

\begin{abstract}

We study the closed group of homeomorphisms
of the boundary of real hyperbolic space generated by a cocompact
Kleinian group $G_1$ and a quasiconformal conjugate $h^{-1}G_2 h$
of a cocompact group $G_2$. We show that if the conjugacy $h$
is not conformal then this group contains a non-trivial one parameter
subgroup. This leads to rigidity results; for example,
Mostow rigidity is an immediate consequence. We are also
able to prove a relative version of Mostow rigidity, called
pattern rigidity. For a cocompact group $G$, by a $G$-invariant pattern
we mean a $G$-invariant collection of closed proper
subsets of the boundary of hyperbolic space which is discrete
in the space of compact subsets minus singletons. Such a pattern arises
for example as the collection of translates of limit sets of
finitely many infinite index quasiconvex subgroups of $G$. We
prove that (in dimension at least three)
for $G_1, G_2$ cocompact Kleinian groups, any quasiconformal map
pairing a $G_1$-invariant pattern to a $G_2$-invariant pattern
must be conformal. This generalizes a previous result of Schwartz
who proved rigidity in the case of limit sets of cyclic subgroups,
and Biswas-Mj \cite{bimj} who proved rigidity for Poincare Duality subgroups.

\smallskip

\begin{center}


\end{center}

\end{abstract}

\overfullrule=0pt

\tableofcontents


\section{Introduction}
The purpose of this paper is to introduce a technique of proving rigidity by constructing
a one-parameter family of homeomorphisms. Our starting data will usually
 be a cocompact Kleinian
group $G$ (i.e. a discrete cocompact group of isometries of $SO(n,1)$) and a generic
quasiconformal map $\phi$ not in $G$. It will be shown that
the closed indiscrete group
$\overline{\langle G, \phi \rangle}$
topologically generated by $G$ and $\phi$ typically contains a flow, i.e.
a one-parameter family of homeomorphisms. This can be thought of as a
`weak Montgomerry-Zippin type theorem for quasiconformal maps'.

\subsection{Statement of results}

Boundaries of Gromov hyperbolic groups are well known examples of
self-similar geometric objects. The self-similarity is a
consequence of the group action on the boundary: arbitrarily small
neighbourhoods of any point can be taken to sets of fixed
diameter. One can "zoom-in" using the group to translate
infinitesimal information into global information leading to
strong rigidity results. This idea goes back to Mostow
\cite{mostow-ihes} in the proof of his celebrated rigidity theorem
for cocompact Kleinian groups in dimension $N \geq 3$. Zooming-in
at a point of differentiability of a quasiconformal conjugacy
between two such groups leads to a linear conjugacy. Similarly in
the theorems of Sullivan \cite{sullivan-rigidity} and Tukia
\cite{tukia} on uniformly quasiconformal groups, and
generalizations of Mostow Rigidity due to McMullen
\cite{ctm-renorm}, zooming-in on a measurable invariant ellipse
field at a point of density of the limit set leads to a constant
invariant ellipse field. Zooming-in is also an essential part of
Schwartz pattern rigidity theorem for symmetric patterns of
geodesics in rank one symmetric spaces \cite{schwarz-inv} and
pattern rigidity for certain quasiconvex Duality and Poincare
duality subgroups of uniform lattices in  rank one symmetric
spaces due to Biswas-Mj \cite{bimj}.\footnote{While Biswas-Mj
\cite{bimj} prove pattern rigidity for quasiconvex duality and PD
subgroups of uniform lattices in  rank one symmetric spaces, this
paper deals with a larger class of subgroups (including {\em all}
quasiconvex subgroups having non-empty domain of discontinuity) of
uniform lattices in {\em real hyperbolic space}. Thus, neither is
subsumed in the other. Besides, the techniques of \cite{bimj} and
the present paper differ widely. While the former relies on
fixed-point theory for homology manifolds, this paper relies on
precise analytical estimates.}

In the present article, we investigate quantitatively
 the consequences of
zooming-in at fixed points of smooth maps. It turns out that
discrete local dynamics translates into continuous global
dynamics. To be precise, for $N \geq 3$, let Homeo($\partial
\mathbb{H}^N$) denote the group of homeomorphisms of the boundary
of $N$-dimensional hyperbolic space equipped with the
uniform topology, $G$ a cocompact Kleinian
group and $f$ a homeomorphism of the boundary. Considering the
upper half-space model of hyperbolic space with boundary
$\mathbb{R}^N \cup \{\infty\}$, when $f$ has a fixed point at
$\infty$ and is tangent to the identity at $\infty$, we have:

\medskip

\begin{theorem} \label{tangent identity}{ If $f(w) = w + \Phi(w) + o(1)$ in a
neighbourhood of $w = \infty$, where $\Phi$ is homogeneous of
degree zero, not identically zero, and satisfies $\Phi(w + v) =
\Phi(w) + O(||v||/||w||)$, then $\overline{<G,f>}$ contains a
non-trivial one-parameter subgroup $(f_t)_{t \in \mathbb{R}}$,
conformally conjugate to a flow of Euclidean translations.}
\end{theorem}

\medskip

The proof consists of zooming-in at $w = \infty$ through a thin
angular sector centered around a ray, where $f$ looks like a translation;
the more one zooms in the smaller the translation looks and
in the limit iterating infinitesimally small translations
one obtains a continuous flow of translations.

\medskip

In the case of fixed points where the derivative is not tangent to
the identity we have the following:

\medskip

\begin{theorem} \label{conformalisable}{ If $f$ is a $C^2$ diffeomorphism
with a fixed point $x_0$ which is not a fixed point of any element
of $G$, and $Df(x_0)$ is conjugate to a conformal linear map
$\lambda O$ with $O$ orthogonal and $\lambda \neq 1$, then
$\overline{<G,f>}$ contains a non-trivial one-parameter subgroup
$(f_t)_{t \in \mathbb{R}}$, conformally conjugate to a flow of
affine linear maps.}
\end{theorem}

\medskip

As with Theorem \ref{tangent identity}, the key to the proof of
the above Theorem is zooming-in with sufficient precision: one
zooms-in near the fixed point of $f$ by group elements with near
by fixed points, but not too large magnification. Then the
commutators of $f$ with these group elements look like real affine
maps close to the identity, and in the limit one iterates such
maps infinitely close to the identity to get a continuous flow of
affine maps.

\medskip

Now consider two cocompact Kleinian groups $G_1, G_2$ and
$h$ a homeomorphism of the boundary (in applications
$h$ will pair $G_i$-invariant structures). By the poles
of a group we mean the set of fixed points of its elements.
This is slightly different from Gromov's  \cite{gromov-hypgps}
definition of poles.
We say $h$ is pole-preserving if it takes poles of $G_1$
to poles of $G_2$. Recall that each group elements has exactly
two fixed points, one attracting, one repelling; if in addition
$h$ takes each such pair of poles of $G_1$ to a pair of
poles of $G_2$ we say $h$ is a pole-pairing map.
The previous Theorem leads naturally to:

\medskip

\begin{theorem} \label{pole-pairing}{ Let $h$ be a $C^2$ diffeomorphism. Then:

\smallskip

\noindent (1) If $\overline{<G_1, h^{-1}G_2 h>}$ does not
contain a non-trivial one-parameter subgroup, then $h$ is
pole-preserving, and the length spectra of $G_1, G_2$
are commensurable.

\smallskip

\noindent (2) If in addition $h$ is linear then $h$ is
also pole-pairing.}
\end{theorem}

\medskip

Statement (1) above is a fairly easy consequence of Theorem
\ref{conformalisable}, whereas statement (2) requires some
computations to see that if $h$ is not pole-pairing then one
obtains a map satisfying the hypotheses of Theorem \ref{tangent
identity}.

\medskip

We apply the above results to the study of the closed subgroup
$\hat{G}$ of homeomorphisms generated by a cocompact Kleinian group $G_1$
and a quasi-conformal conjugate $h^{-1} G_2 h$ of a cocompact
Kleinian group $G_2$. Our main theorem is the following:

\medskip

\begin{theorem}\label{main theorem}{ If $h$ is a quasi-conformal map
which is not conformal, then the group $\hat{G} = \overline{<G_1,
h^{-1}G_2 h>}$ contains a non-trivial one-parameter subgroup (and
in particular is not discrete).}
\end{theorem}

\medskip

The proof proceeds by first performing a preliminary zoom-in to
upgrade the quasi-conformal map $h$ to a non-conformal linear map
$A$, and using Theorem \ref{pole-pairing} to see that $A$ is
pole-pairing. One then picks non-linear conformal maps $g_i \in
G_i$ whose poles correspond under $A$, and conjugates the groups
$G_i$ by conformal maps sending the poles of the $g_i$'s to $0,
\infty$ to get groups $G_i'$, conformal linear maps $g_i' \in
G_i'$, and a group $\hat{G}' = \overline{<G_1', \mu^{-1}G_2'
\mu>}$ where $\mu$ is a smooth pole-pairing map fixing $0,\infty$
given by pre and post composing $A$ with conformal maps. The key
observation is that the non-conformality of $A$ implies that $\mu$
is not linear, a consequence of the fact that the only linear maps
conjugate to linear maps by inversions in spheres around the
origin are conformal linear maps. So zooming-in in at the fixed
point $0$ of $\mu$ by $g_1'$ and zooming-out by $g_2'$ one gets a
sequence of non-linear pole-pairing maps $h_n$ converging to a
linear pole-pairing map $B$, so in particular a sequence
containing infinitely many distinct maps. For any $g \in G_2'$
with poles $a,b$ the conjugates $f_n = h_n^{-1} g h_n$ belong to
$\hat{G}'$ and have fixed points $a_n = h_n^{-1}(a), b_n =
h_n^{-1}(b)$ converging to $a' = B^{-1}a, b' = B^{-1}b$. The pairs
$(a_n,b_n)$ are poles of some $g_n \in G_1'$ and zooming-in on
$f_n$ by $g_n$ one obtains maps $F_n \in \hat{G}'$ which are
conformal conjugates of linear maps, with fixed points $a_n, b_n$,
converging to $F \in \hat{G}'$ a conformal conjugate of a linear
map, with fixed points $a',b'$. Using a "Scattering Lemma" due to
Schwartz \cite{schwarz-inv} and density of poles, there must be a
$g \in G_2'$ such that $a_n \neq a', b_n \neq b'$ for infinitely
many $n$. Then as in the proof of Theorem \ref{conformalisable},
one considers the commutators of the maps $F_n$ and $F$ to get in
the limit maps looking like affine maps infinitely close to the
identity, which one iterates to get a continuous flow.

\medskip

Note Mostow Rigidity (where $G_1 = h^{-1} G_2 h$ so $\hat{G} =
G_1$ is discrete) is an immediate consequence. Taking $G_1 = G_2 =
G$ say in Theorem \ref{main theorem} gives

\medskip

\begin{cor} \label{one qc map}{ If $h$ is a quasi-conformal map which is not
conformal, then $\overline{<G, h>}$ contains a non-trivial
one-parameter subgroup.}
\end{cor}

\medskip

The main application of Theorem \ref{main theorem} is to the
problem of {\it pattern rigidity}. This problem was motivated in
part by work of Mosher-Sageev-Whyte \cite{msw2} on quasi-isometric
rigidity for fundamental groups of graphs of groups with vertex
groups cocompact Kleinian groups and edge groups quasiconvex
subgroups of the adjacent vertex groups. In \cite{msw2} it is
shown that a quasi-isometry between two such groups leads to
pairings of {\it symmetric patterns of limit sets}, by which we
mean the following:

Given $G$ a uniform lattice in a rank one symmetric space, a $G$-symmetric pattern of
limit sets $\mathcal{J}$ is a $G$-invariant collection of translates of
limit sets of finitely many
infinite index quasiconvex subgroups $H_1,\dots, H_n$.
The quasiconvexity hypothesis ensures that
such a collection is discrete in the Hausdorff
topology on the space of compact subsets of the boundary minus singletons.

Schwartz' rigidity theorem for symmetric patterns of geodesics \cite{schwarz-inv}
in a rank one symmetric space
can be formulated in this context as rigidity for symmetric patterns
of limit sets when the subgroups $H_1, \dots, H_n$ are cyclic. Biswas-Mj \cite{bimj}
generalize this to the case of subgroups which are either
codimension one duality groups or odd-dimensional Poincare Duality
groups or subgroups containing such subgroups as free factors.

 We obtain as an easy corollary of Theorem \ref{main theorem} the general case for uniform lattices
in $SO(n,1)$:

\medskip

\begin{cor} \label{pattern rigidity}{ For $N \geq 3$, suppose
$\mathcal{J}_i$ (for $i = 1,2$) are $G_i$-invariant collections of
closed subsets of $\partial \mathbb{H}^N$ which are discrete in
the Hausdorff topology on the space of compact subsets of the boundary minus singletons.
If $h(\mathcal{J}_1) = \mathcal{J}_2$ for a quasi-conformal homeomorphism $h$, then $h$ must be conformal.}
\end{cor}

\medskip

%

In particular we get pattern rigidity for symmetric patterns of
limit sets. It also follows (using the main result of \cite{mahan-agt})
that if $\mathcal{J}_i, i=1,2$ are
symmetric patterns of closed convex (or quasiconvex) sets in hyperbolic space
$\mathbb{H}^N$, then any uniformly proper bijection between
$\mathcal{J}_1$ and $\mathcal{J}_2$ is induced by a hyperbolic
isometry (for the definition of uniformly proper we refer to \cite{schwarz-inv}, \cite{mahan-agt}).
Combining pattern rigidity with the main theorem of
Mosher-Sageev-Whyte \cite{msw2} (to which we refer for the terminology)
we get the following quasi-isometric rigidity theorem:

\medskip

\begin{cor} \label{qirig} { Let $\mathcal{G}$ be a
finite, irreducible graph of groups with associated Basse-Serre
tree $T$ of spaces such that no depth zero raft of $T$ is a line.
Further suppose that the vertex groups are fundamental groups of
compact hyperbolic manifolds of dimension $N \geq 3$ and edge
groups are infinite index quasiconvex subgroups of the adjacent
vertex groups.

\smallskip

If $H$ is a finitely-generated group quasi-isometric to $G =
\pi_1(\mathcal{G})$ then $H$ splits as a graph of groups
$\mathcal{G'}$ whose depth zero vertex groups are commensurable to
those of $\mathcal{G}$ and whose edge groups and positive depth
vertex groups are quasi-isometric to those of $\mathcal{G}$.}
\end{cor}

\bigskip

%

\subsection{Acknowledgements}
I thank Mahan Mj for numerous helpful
discussions and for encouragement to venture into unfamiliar
territory. I am grateful to Marc Bourdon and Herve Pajot for beautiful lectures on
quasi-isometric rigidity and geometric function theory.

\section{Zooming-in: Quantitative Estimates}

Throughout this section and the next $G$ will denote a cocompact Kleinian
group. As mentioned before, we shall refer to a discrete subgroup of $SO(n,1)$
as a Kleinian group. We first recall an elementary result on zooming-in by
a group element at a fixed point of a differentiable map in 2.1, and
then develop estimates for zooming-in by a sequence of group elements
near a fixed point in 2.2 and 2.3. We work throughout with the upper-half space model
$\mathbb{H}^N = \{ w = (x,t) | x \in \mathbb{R}^{N-1}, t > 0
\}$ with boundary $\partial \mathbb{H}^N = \mathbb{R}^{N-1} \cup
\{\infty\}$.

\subsection{Zooming-in at fixed points}

Recall that the boundary map of a hyperbolic translation is
a conformal map with fixed points $\epsilon \neq
M$, given by the conjugate of a similarity
(or conformal linear map) $\lambda O$ ($\lambda \neq 1, O$
orthogonal) by a conformal map sending $0$ to $\epsilon$, $\infty$
to $M$. Zooming-in on quasiconformal maps at fixed points
which are points of differentiability will give us maps of the
following form:

\begin{defn}\label{based linear map}{ Let $\epsilon, M \in \mathbb{R}^{N-1} \cup \{\infty\},
\epsilon \neq M$. By a {\bf linear map based at $\epsilon, M$ with
multiplier $A$} we mean a map $g$ which is the conjugate of a
linear map $A$ by a conformal map $S$ sending $0$ to $\epsilon$,
$\infty$ to $M$ such that $DS(0) = I$ (so that $Dg(\epsilon) =
A$).}
\end{defn}

\medskip

\begin{prop}\label{zoom-in at pole}{ Let $f$ be a homeomorphism of
$\mathbb{R}^{N-1} \cup \{\infty\}$ with a fixed point $\epsilon$
such that $f$ is differentiable at $\epsilon$. Let $T$ be the
boundary map of a hyperbolic translation with fixed points
$\epsilon, M$ such that $||DT(\epsilon)|| > 1$. Then along a
subsequence $T^n f T^{-n}$ converges uniformly on compacts of
$\mathbb{R}^{N-1} \cup \{\infty\} - M$ to a linear map based at
$\epsilon, M$ with multiplier $Df(\epsilon)$.}
\end{prop}

\medskip

\noindent {\bf Proof:} Let $DT(\epsilon) = \lambda O, \lambda > 1$.
If $\epsilon = 0$ then along a subsequence such that
$O^n \to I$,
\begin{align*}
\lambda^n O^n f(O^{-n}x/\lambda^n) & = \lambda^n O^n
(Df(0)(O^{-n}x/\lambda^n) + o(x/\lambda^n)) \\ & = O^n Df(0) O^{-n}(x)
+ o(x) \to Df(0)(x) \\
\end{align*}
uniformly on compacts of $\mathbb{R}^{N-1}$. The general case
follows on conjugating $f, T$ by a conformal map sending
$\epsilon$ to $0$, $M$ to $\infty$. $\diamond$

\subsection{A quantitative estimate for density of poles}

It is well known that
poles of $G$ (i.e. fixed points of group elements) are dense in
the space of pairs of points on the boundary $\partial
\mathbb{H}^N$. We will need the following effective bound for the error
in approximating a given point by poles in terms of the translation
lengths of the corresponding group elements (taking the point
to be the origin in $\mathbb{R}^{N-1}$):

\medskip

\begin{prop} \label{pole density} { There are sequences $0 < t_n
\leq 1, g_n \in G$ such that the poles $\epsilon_n, M_n$ and
translation lengths $l_n$ of the $g_n$'s satisfy $l_n \to +\infty,
||\epsilon_n|| \ll t_n e^{-l_n}, ||M_n|| \gg t_n$.}
\end{prop}

\medskip

\noindent (We write $a_n \ll b_n$ for $a_n = o(b_n)$ and
$a_n \lesssim b_n$ for $a_n = O(b_n)$. We denote the geodesic segment between
points $u,v \in \mathbb{H}^N$ by $[u,v]$.)

\medskip

\noindent {\bf Proof:} Let $\pi : \mathbb{H}^N \to M = \mathbb{H}^N / G$ denote the
universal covering map. By compactness of the unit tangent bundle of $M$,
the unit tangent vectors to the semi-infinite geodesic
$\pi(\gamma = \{(0,e^{-l})| l \geq 0 \})$ must have an accumulation point as $l \to +\infty$.
It follows that we can find sequences $0 < t_n \leq 1, l_n' \to +\infty$
and $g_n \in G$
such that $\delta_n = d(g_n(0,t_n), (0,t_n e^{-l_n'})) \to 0$, and the angle
$\theta_n$ between the vertical direction (pointing downwards) and $dg_n((0,t_n))(\overrightarrow{u_n})$
tends to $0$, where $\overrightarrow{u_n}$ is the unit
tangent vector to $\gamma$ at the point $(0,t_n)$. Fix $\epsilon >
0$ such that $1/2 < \sin(\theta)/\theta, \tan(\theta)/\theta < 2$
for $0 < \theta < \epsilon$. Without loss of generality
we may assume that $\theta_n < \epsilon/2$ and $\lambda_n'= e^{-l_n'} < 1/2$ for all $n$.
Let $\epsilon_n, M_n$ be the attracting and repelling poles respectively of
$g_n$, and $l_n$ the translation length of $g_n$.

\medskip

Let $w_{n,1} := (0, t_n), w_{n,2} := (0, t_n e^{-l_n'})$ and, for $j
\geq 1$, let $w_{n,2j+1} = (x_{n, 2j+1}, t_{n,2j+1}) :=
g_n^j(w_{n,1})$, $w_{n,2j+2} = (x_{n, 2j+2}, t_{n,2j+2}) :=
g_n^j(w_{n,2})$. For $j \geq 1$ let $\theta_{n,j}$ be the angle between
the vertical direction and the geodesic segment $[w_{n,2j+1}, w_{n,2j+2}]$
at $w_{n,2j+1}$ (so $\theta_{n,1} = \theta_n$), $\theta_{n,j}'$
the angle between the tangent vector to $[w_{n,2j+1}, w_{n,2j+2}]$ at
$w_{n,2j+2}$ and the tangent vector to $[w_{n,2j+3}, w_{n,2j+4}]$ at
$w_{n,2j+3}$ (considered as elements of $\mathbb{R}^N$), and
$\tau_{n,j}$ the angle between the vertical
direction and $[w_{n,2j+1}, w_{n,2j+2}]$ at $w_{n,2j+2}$. Let
$\theta_n'$ be the angle at $w_{n,2j+3}$ between the tangent to
$[w_{n,2j+3}, w_{n,2j+4}]$ and the tangent to
$[w_{n,2j+1},w_{n,2j+2}]$ at $w_{n,2j+2}$ parallel transported to $w_{n,2j+3}$
along $[w_{n,2j+2},w_{n,2j+3}]$ (note this is independent of $j$).
Since $d(w_{n,2j+2}, w_{n,2j+3}) = \delta_n \to 0$,
$|\theta_{n,j}' - \theta_n'|, |\theta_n' - \theta_n| \leq \eta_n/2$ for some sequence
$\eta_n \to 0$, so $|\theta_{n,j}' - \theta_n| \leq \eta_n$.
We will prove
by induction on $j$ that, for $n$ sufficiently large, the following hold
for all $j \geq 1$:

\medskip

\noindent (1) $t_{n,2j} \leq 2^{j-1}{\lambda_n'}^{(j+1)/2}t_n$

\medskip

\noindent (2) $\theta_{n,j} \leq 2 [(\theta_n + \eta_n)(1 + 4 {\lambda_n'}^{1/2}
+ (4{\lambda_n'}^{1/2})^2 + \dots + (4{\lambda_n'}^{1/2})^{j-2}) + 2^{2j-3}{\lambda_n'}^{(j-1)/2}\theta_n]$.

\medskip

The inequalities clearly hold for $j=1$ since $t_{n,2} =
\lambda_n' t_n, \theta_{n,1} = \theta_n$. Assume they hold for
some $j \geq 1$. Since $d(w_{n,2j}, w_{n,2j+1}) =
d(w_{n,2},w_{n,3}) = \delta_n \to 0$ as $n \to \infty$, clearly
$t_{n,2j+1} \leq 2 t_{n,2j}$ for $n$ large enough. Integrating the
hyperbolic metric along the straight line segment joining
$w_{n,2j+1}$ to $w_{n,2j+2}$ gives the estimate
$$
l_n' = d(w_{n,2j+1}, w_{n,2j+2}) \leq
\log\left(\frac{t_{n,2j+1}}{t_{n,2j+2}}\right) (1 +
\tan(\theta_{n,j}))
$$
For $n$ large enough the RHS of (2) will be small enough for us to
assume $\tan(\theta_{n,j}) < 1$, and thus
\begin{align*}
t_{n,2j+2} & \leq t_{n,2j+1} e^{-l_n'/(1 + \tan(\theta_{n,j}))} \\
           & \leq 2 t_{n,2j} {\lambda_n'}^{1/2} \\
           & \leq 2
           \left(2^{j-1}{\lambda_n'}^{(j+1)/2}t_n\right){\lambda_n'}^{1/2}
           = 2^j{\lambda_n'}^{(j+2)/2}t_n \\
\end{align*}

This proves (1). For (2), we consider Figure 1 below,

\medskip

{\hfill {\centerline {\psfig {figure=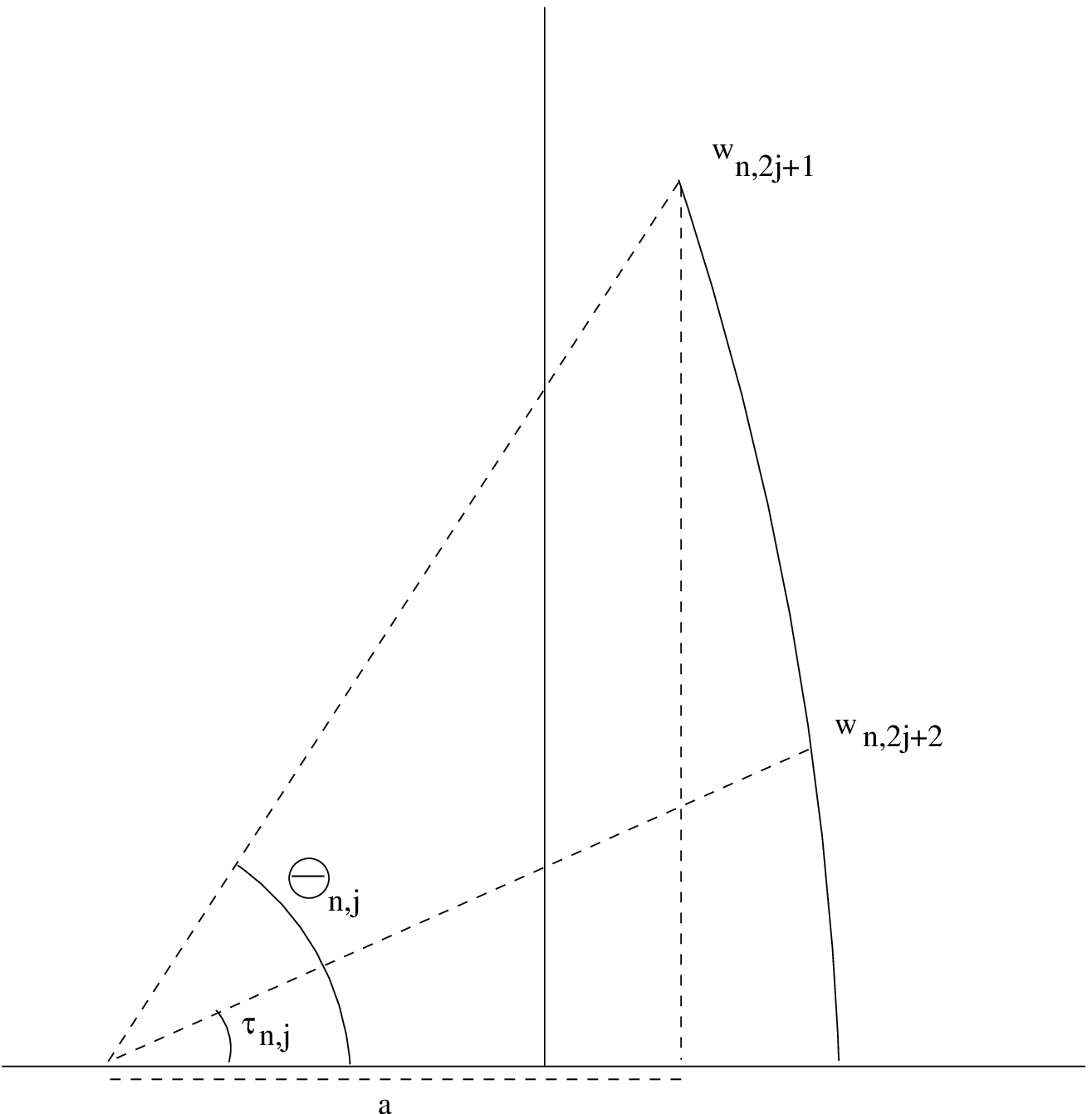,height=8cm}}}}

{\centerline {\bf Figure 1}}

\medskip

from which
we see that $\tan(\tau_{n,j}) \leq t_{n,2j+2}/a,
\tan(\theta_{n,j}) = t_{n,2j+1}/a$, so
$\tan(\tau_{n,j})/\tan(\theta_{n,j}) \leq t_{n,2j+2}/t_{n,2j+1} \leq {\lambda_n'}^{1/2}$;
for $n$ large enough, the angles $\tau_{n,j}, \theta_{n,j}$ are thus small
enough to be comparable with their tangents, hence
$\tau_{n,j}/\theta_{n,j} \leq 2 {\lambda_n'}^{1/2}$. For unit
vectors $u, v$ in $\mathbb{R}^N$ at an angle $\theta$ to each other,
the Euclidean norm $||u - v||$ is equal to $2 \sin(\theta/2)$, so
by the triangle inequality we have $\sin(\theta_{n,j+1}/2) \leq
\sin(\theta_{n,j}'/2) + \sin(\tau_{n,j}/2)$, which (for $n$ large enough
so all the angles above are small enough) implies
\begin{align*}
\theta_{n,j+1} & \leq 2(\theta_{n,j}' + \tau_{n,j}) \\
               & \leq 2(\theta_n + \eta_n+
               2{\lambda_n'}^{1/2}\theta_{n,j}) \\
& \leq 2[(\theta_n + \eta_n)+2{\lambda_n'}^{1/2} 2 [(\theta_n + \eta_n)(1 + 4 {\lambda_n'}^{1/2}
+ \dots + (4{\lambda_n'}^{1/2})^{j-2}) + 2^{2j-3}{\lambda_n'}^{(j-1)/2}\theta_n]]\\
& \leq 2[(\theta_n + \eta_n)(1+4{\lambda_n'}^{1/2} + \dots +
(4{\lambda_n'}^{1/2})^{j-1}) + 2^{2j-1}{\lambda_n'}^{j/2}\theta_n]\\
\end{align*}
This proves (2), and it follows that for all $n$ large enough, we will have
$\tan(\theta_{n,j}) \leq C (\theta_n+\eta_n)$ for all $j$, for some constant $C$.
Considering the tangent line to $[w_{n,2j+1},w_{n,2j+2}]$ at
$w_{n,2j+1}$ gives $||x_{n,2j+2} - x_{n,2j+1}|| \leq
\tan(\theta_{n,j}) t_{n,2j+1}$. We also have $||x_{n,2j+1} - x_{n,2j}|| \leq
C' \delta_n t_{n,2j}$ for some constant $C'$, hence
\begin{align*}
||\epsilon_n|| & = ||\lim_{j \to +\infty} g_n^j(w_{n,1})|| \\
           & = ||\lim_{j \to +\infty} x_{n,j}|| \\
           & \leq \sum_{j = 1}^{\infty} ||x_{n,2j+1} - x_{n,2j}||
           + ||x_{n,2j+2} - x_{n,2j+1}|| \\
           & \leq \sum_{j = 1}^{\infty} C' \delta_n 2^{j-1}{\lambda_n'}^{(j+1)/2}t_n
           + 2C(\theta_n + \eta_n)
           2^{j-1}{\lambda_n'}^{(j+1)/2}t_n \\
           & \leq C'' \beta_n \lambda_n' t_n \\
\end{align*}
for some constant $C''$ and some sequence $\beta_n \to 0$. The
translation length $l_n$ of $g_n$ satisfies $l_n \leq l_n' +
\delta_n$, hence $e^{-l_n} \gtrsim e^{-l_n'} = \lambda_n'$, thus
$||\epsilon_n|| \ll t_n e^{-l_n}$.

\medskip

For the estimate on the other
pole $M_n$, we conjugate by inversion in the unit sphere and apply
an argument similar to the above to get $1/||M_n|| \ll 1/t_n$, hence
$||M_n|| \gg t_n$. $\diamond$

\bigskip

\subsection{Almost affine maps}

\bigskip

We suppose as above that $0 \in \mathbb{R}^{N-1}$
is not a fixed point of any element of $G$, and consider the
sequence $g_n$ given by the
previous Proposition. We now adjust the poles, conjugating by
appropriately chosen
dilations to move $M_n$ closer to $\infty$ and $\epsilon_n$
away from $0$, in order to get a sequence of maps which look like
affine maps up to first order:

\medskip

\begin{prop} \label{almost affine maps} { There exist dilations $T_n : x
\mapsto x/t_n', t_n'>0$, such that the following hold for the maps
$g_n' = T_n g_n T_n^{-1}$:

\smallskip

\noindent (1) $1/||M_n'|| \ll ||\epsilon_n'||/\lambda_n \ll 1, \ t_n' \ll
||\epsilon_n'||/\lambda_n, \ ||\epsilon_n'|| \gg \lambda_n^2$

\smallskip

\noindent (where
$\epsilon_n' = \epsilon_n/t_n', M_n' = M_n/t_n'$ are the poles
of $g_n'$, and $\lambda_n = e^{-l_n}$).

\smallskip

\noindent (2) $g_n'(x) = \lambda_n O_n(x) +
\epsilon_n' + o(\epsilon_n')$

\smallskip

\noindent (where $O_n$ is the orthogonal
part of the derivative of $g_n$ at $\epsilon_n$, and the estimate
holds uniformly on compacts).
}
\end{prop}

\medskip

\noindent {\bf Proof:} The inequalities in (1) are equivalent
to the inequalities
$$
\frac{||\epsilon_n||}{\lambda_n} \ll t_n' \ll \min\left(
\left(\frac{||\epsilon_n|| ||M_n||}{\lambda_n}\right)^{1/2}, \left(
\frac{||\epsilon_n||}{\lambda_n}\right)^{1/2},
\frac{||\epsilon_n||}{\lambda_n^2}\right)
$$
for $t_n'$; by the estimates of the previous Proposition it is
easily seen that the $LHS \ll RHS$ above, hence $t_n'$ can be
chosen appropriately to satisfy the above inequalities.

\medskip

For the estimate (2), we use the following explicit formula
for $g_n'$,
$$
g_n'(x) = \frac{\lambda_n O'_n \xi - (M_n' - \epsilon_n')}{||\lambda_n O'_n
\xi - (M_n' - \epsilon_n')||^2} ||M_n' - \epsilon_n'||^2 + M_n'
$$
where $O'_n = \rho_n O_n \rho_n$, $\rho_n$ being the
reflection in the hyperplane normal to $M_n'$, and
$$
\xi = \frac{x - M_n'}{||x - M_n'||^2} ||M_n' - \epsilon_n'||^2 +
M_n' - \epsilon_n'
$$
Using the inequalities in (1), the formula $||a - b||^2 = ||a||^2 -
2<a,b> + ||b||^2$ and the geometric series gives
$$
\frac{||M_n' - \epsilon_n'||^2}{||x - M_n'||^2} = 1 + 2\frac{<x -
\epsilon_n', M_n'>}{||M_n'||^2}+ O(1/||M_n'||^2)
$$
and hence, by substituting in the formula for $\xi$,
\begin{align*}
\xi & = (x - \epsilon_n') - 2\frac{<x -
\epsilon_n', M_n'>}{||M_n'||^2}M_n' + O(1/||M_n') \\
    & = \rho_n(x - \epsilon_n') + O(1/||M_n'). \\
\end{align*}
Similarly,
$$
\frac{||M_n' - \epsilon_n'||^2}{||\lambda O'_n \xi - (M_n' -
\epsilon_n')||^2} = 1 + 2\frac{<\lambda_n O'_n \xi, M_n'>}{||M_n'||^2}
+ O(1/||M_n'||^2)
$$
so substituting in the formula for $g_n'$ and using
(1) gives,
\begin{align*}
g_n'(x) & = \lambda_n O'_n \xi - 2\frac{<\lambda_n O'_n \xi,
M_n'>}{||M_n'||^2}M_n' + \epsilon_n' + o(\epsilon_n') \\
        & = \rho_n(\lambda_n O'_n \xi) + \epsilon_n' + o(\epsilon_n')
        \\
        & = \rho_n(\lambda_n O'_n \rho_n(x - \epsilon_n')) +
        \epsilon_n' + o(\epsilon_n') \\
        & = \lambda_n O_n (x) + \epsilon_n' +
        o(\epsilon_n')
\qquad \qquad \qquad \qquad \qquad \diamond\\
\end{align*}

We now estimate what we see when we zoom-in
using these ''almost
affine'' maps $g_n'$ on a $C^2$ diffeomorphism $f$ with a fixed point at
the origin (which we continue to assume is not a
pole of $G$):

\medskip

\begin{prop} \label{almost affine zoom}{ Let $B = Df(0)$ and $f_n = T_n f
T_n^{-1}$. Then:

\smallskip

\noindent (1) $g_n'^{-1}f_n g_n'(x)= O_n^{-1}BO_n(x) + O_n^{-1}(B -
I)(\epsilon_n'/\lambda_n) + o(||\epsilon_n'||/\lambda_n)$

\smallskip

\noindent (2) $f_n^{-1}g_n'^{-1}f_n g_n'(x)= B^{-1}O_n^{-1}BO_n(x) + B^{-1}O_n^{-1}(B -
I)(\epsilon_n'/\lambda_n) + o(||\epsilon_n'||/\lambda_n)$
}
\end{prop}

\medskip
\medskip

\noindent {\bf Proof:} Using the Taylor expansion $f(h) = Bh +
O(||h||^2)$ near $0$, we have
\begin{align*}
f_n(g_n'(x)) & = \frac{1}{t_n'}f(t_n'g_n'(x)) \\
             & = \frac{1}{t_n'}\left(B(t_n'g_n'(x)) +
             O(t_n'^2||g_n'(x)||^2)\right) \\
             & = B(g_n'(x)) + O(t_n'||g_n'(x)||^2) \\
             & = \lambda_n B O_n x + B \epsilon_n' +
             o(\epsilon_n') \\
\end{align*}
(using the estimates (1),(2) of the previous
Proposition). In terms of the variables
$$
x' = f_n(g_n'(x)), \ \xi = (x' - M_n')\frac{||M_n' -
\epsilon_n'||^2}{||M_n' - x'||^2} + (M_n' - \epsilon_n'),
$$
we have
$$
g_n'^{-1} f_n g_n' (x) = (\lambda_n^{-1} O_n'^{-1}\xi - (M_n' -
\epsilon_n'))\frac{||M_n' -
\epsilon_n'||^2}{||\lambda_n^{-1}O_n'^{-1}\xi - (M_n' -
\epsilon_n')||^2} + M_n'
$$
Using the fact that $x' = O(\lambda_n)$ we get
$$
\frac{||M_n' -
\epsilon_n'||^2}{||M_n' - x'||^2} = 1 + 2 \frac{<x' - \epsilon_n',
M_n'>}{||M_n'||^2} + O\left(\frac{\lambda_n}{||M_n'||^2}\right)
$$
which leads to the estimates
\begin{align*}
\xi                        & = \rho_n(x' - \epsilon_n') + O(\lambda_n/||M_n'||), \\
\lambda_n^{-1}O_n'^{-1} \xi & = O_n'^{-1} \rho_n((x' - \epsilon_n')/\lambda_n)+O(1/||M_n'||) \\
                           & = O_n'^{-1} \rho_n B O_n' x + O_n'^{-1}
                           \rho_n (B - I) (\epsilon_n'/\lambda_n) +
                           o(||\epsilon_n'||/\lambda_n) \\
\end{align*}
(using $1/||M_n'|| \ll ||\epsilon_n'||/\lambda_n)$).
Similar calculations then give
\begin{align*}
g_n'^{-1}f_n g_n'(x) & = \rho_n O_n'^{-1} \rho_n B O_n x + \rho_n O_n'^{-1}
                           \rho_n (B - I) (\epsilon_n'/\lambda_n) +
                           o(||\epsilon_n'||/\lambda_n) \\
                     & = O_n^{-1}B O_n x + O_n^{-1}
                           (B - I) (\epsilon_n'/\lambda_n) +
                           o(||\epsilon_n'||/\lambda_n) \\
\end{align*}
Finally, writing $\eta = g_n'^{-1}f_n g_n'(x)$, we have
\begin{align*}
f_n^{-1}g_n'^{-1}f_n g_n'(x) = f_n^{-1}(\eta) & = \frac{1}{t_n'}(B^{-1}(t_n'\eta) + O(t_n'^2||\eta||^2)) \\
                                              & = B^{-1}\eta +
                                              O(t_n') \\
                                              & = B^{-1}O_nBO_n(x) + B^{-1}O_n^{-1}(B -
I)(\epsilon_n'/\lambda_n) + o(||\epsilon_n'||/\lambda_n) \\
\end{align*}
(using $t_n' \ll ||\epsilon_n'||/\lambda_n$). $\diamond$

\bigskip

\section{Generating flows}

\subsection{Euler's formula and affine maps}

\bigskip

For $N \geq 2$, we denote by Aff(${\mathbb{R}}^N$) the Lie group
of real affine maps of ${\mathbb{R}}^N$ (given by maps of the form
$x \mapsto Bx + b, B \in GL_N({\mathbb{R}}), b \in
{\mathbb{R}}^N$), and by {\it aff} (${\mathbb{R}}^N$) its Lie
algebra, given by real affine vector fields $x \mapsto A(x) = Bx +
b, B \in M_N({\mathbb{R}}), b \in {\mathbb{R}}^N$. We equip {\it
aff} (${\mathbb{R}}^N$) with the norm $||A|| := Max(||B||,
||b||)$, where $A(x) = Bx + b$. We write $e^A$ for the exponential
of $A \in $ {\it aff}$({\mathbb{R}}^N)$. The exponential is a
local diffeomorphism near $0$, and the local inverse satisfies $$
\frac{||\log(id + A)||}{||A||} \to 1 $$ when $A \in $ {\it
aff}$({\mathbb{R}}^N)$ tends to $0$.

%

\medskip

Recall Euler's
formula,
$$
\left( 1 + \frac{x}{n} + o\left(\frac{1}{n}\right) \right)^n \to
e^x.
$$
We will need the following version of this formula
when $x \in $ {\it aff}$({\mathbb{R}}^N)$ :

\medskip

\begin{prop} \label{Euler's formula}{ Let $(f_n)$ be a sequence of
maps from ${\mathbb{R}}^N$ to itself such that $f_n = id + A_n +
E_n$, where $A_n \in $ {\it aff} (${\mathbb{R}}^n$), $A_n \neq 0,
||A_n|| \to 0$, and $||E_n(x)||/||A_n|| \to 0$ uniformly on
compacts. Then there is a subsequence of $(f_n)$ and an $A \in $
{\it aff}(${\mathbb{R}}^N$), $A \neq 0$, such that for all $t >
0$, there is a sequence of integers $m_n = m_n(t)$ such that
${f_n}^{m_n} \to e^{tA}$ uniformly on compacts along the
subsequence.}
\end{prop}

\medskip

\noindent {\bf Proof:} Let $g_n = id + A_n$, then
for $n$ sufficiently large $g_n \in$ Aff$({\mathbb{R}}^N)$ and we can write
$g_n = e^{A_n'}$ for some $A_n' \in $ {\it aff}$({\mathbb{R}}^N)$. Let
$t_n = ||A_n'||$; by the formula for the logarithm given above,
$||A_n||/2 \leq t_n \leq 2||A_n||$ for $n$ large enough. The unit
ball in {\it aff}(${\mathbb{R}}^N$) being compact, we can choose a
subsequence such that $A_n'/t_n$ converges to some $A$ along the
subsequence.

\medskip

Given $t > 0$, let $m_n = [t/t_n]$ be the integer part of
$t/t_n$. Then along the chosen subsequence,
$$
g_n^{m_n} = e^{t_n m_n (A_n'/t_n)} \to e^{tA}
$$
uniformly on compacts since $t_n m_n \to t, A_n'/t_n \to A$;
so it suffices to show
that $f_n^{m_n} - g_n^{m_n} \to 0$ uniformly on compacts. We note also that the
elements $g_n^j, n \geq 1, 0 \leq j \leq m_n$, are contained in a compact subset
of Aff$({\mathbb{R}}^N)$ so their derivatives are uniformly bounded on compact subsets
of ${\mathbb{R}}^N$.
Fix a compact subset of ${\mathbb{R}}^N$ and an upper bound $C > 1$ for
these derivatives on the compact.
Given $\epsilon > 0$, by hypothesis on any compact we have $||g_n(x) - f_n(x)|| =
||E_n(x)|| \leq \epsilon ||A_n||$ for $n$ large enough.
Assume for some $1 \leq j \leq m_n - 1$ we have
$||g_n^j(x) - f_n^j(x)|| \leq C j \, \epsilon ||A_n||$. Then by the
Mean Value Theorem and our induction hypothesis we have
\begin{align*}
||g_n^{j+1}(x) - f_n^{j+1}(x)|| & = ||g_n^j(g_n(x)) -
g_n^j(f_n(x)) + g_n^j(f_n(x)) - f_n^j(f_n(x))|| \\
& \leq C ||g_n(x) - f_n(x)|| + C j \, \epsilon ||A_n|| \\
& \leq C \epsilon ||A_n|| + C j \, \epsilon ||A_n|| =
C(j+1)\, \epsilon||A_n|| \\
\end{align*}
so by induction it follows that $||g_n^{m_n}(x) - f_n^{m_n}(x)||
\leq C m_n \, \epsilon ||A_n|| \leq C m_n \, \epsilon (2 t_n) \leq 2Ct\,\epsilon$.
$\diamond$

\subsection{Existence of flows}

Let $G$ be a cocompact Kleinian group.

\medskip

\begin{lemma} \label{cocompact action}{ Let $T_n$  be a sequence of conformal
maps and $f$ a homeomorphism of $\partial \mathbb{H}^N$. Then:

\smallskip

\noindent (1) There is a conformal map $\phi = \phi(G, \{T_n\})$
and $g_n \in G$ such that $g_n^{-1} T_n \to \phi$ along a
subsequence.

\smallskip

\noindent (2) If $h_n \in T_n <G, f> T_n^{-1} \leqslant $Homeo($\partial
\mathbb{H}^N$) converges to $h$ then the closed subgroup
$\overline{<G, f>} \leqslant $Homeo($\partial \mathbb{H}^N$)
generated by $G$ and $f$ contains \, $\phi^{-1} h \phi$. }
\end{lemma}

\smallskip

\noindent {\bf Proof:} Since $G$ acts cocompactly on triples
on the boundary, we can choose $g_n \in G$ such that
the conformal maps $\phi_n = T_n \circ g_n$ uniformly
separate three chosen points
and hence form an equicontinuous family (since cross-ratios
are preserved). So $\phi_n$ converges to a conformal
map $\phi$ along a subsequence. Let $k_n \in <G, f>$
be such that $h_n = T_n k_n T_n^{-1}$. Then the maps
$g_n^{-1} k_n g_n = \phi_n^{-1} h_n \phi_n$
converge to $\phi^{-1} h
\phi$ along the same subsequence. $\diamond$

\medskip

\noindent {\bf Proof of Theorem \ref{tangent identity}:} Take
$w_0$ such that $c = \Phi(w_0) \neq 0$. In thin sectors near
infinity $U_{\epsilon, R} = \{ w = tw_0 + v : t \geq R, ||v|| \leq
\epsilon ||w|| \}$ centered around the ray with direction $w_0$,
the hypotheses on $\Phi$ imply that for large $R$,
\begin{align*}
f(w) & = w + \Phi(w) + o(1) \\
     & = w + \Phi(tw_0 + v) + o(1) \\
     & = w + c + O(\epsilon) + o(1) \\
\end{align*}
Fix positive sequences $\epsilon_n \to 0, R_n \to +\infty$. For $n
\geq 1$ let $T_n(w) = w + a_n w_0$ be a translation in the direction $w_0$
with $a_n > 0$ large enough so that $T_n$ sends the ball of radius
$n$ around the origin into the sector $U_{\epsilon_n,R_n}$, and $S_n(w) = nw$
a scaling factor $n$ (note $S_n$ preserves the sectors $U_{\epsilon_n,R_n})$.
Then, on the ball of radius $n$ around the origin,
\begin{align*}
f_n := T_n^{-1} S_n^{-1} f S_n T_n(w) & = \frac{1}{n}f(n(w+ a_n w_0)) -
a_n w_0 \\
& = \frac{1}{n} (n(w + a_n w_0) + c + O(\epsilon_n) + o(1)) - a_n
w_0 \\
& = w + \frac{c}{n} + \frac{1}{n}(O(\epsilon_n) + o(1)) \\
& = w + \frac{c}{n} + o\left(\frac{c}{n}\right) \\
\end{align*}
from which it follows as in the proof of Proposition \ref{Euler's
formula} that for any $t > 0$, $f_n^{[nt]}(w) \to w + ct$
uniformly on compacts. The Theorem then follows from Lemma
\ref{cocompact action} $\diamond$.

\medskip

\noindent {\bf Proof of Theorem \ref{conformalisable}:} Without
loss of generality we may assume the fixed point of $f$ is at $0$.
Let $B = Df(0)$. Consider the sequences $g_n \in G, T_n, g_n' =
T_n g_n T_n^{-1}$ given by the previous section. Passing to a
subsequence if necessary we may assume $O_n$ converges to some
orthogonal linear map $O$. We consider two cases:

\medskip

\noindent {\it Case 1: Some power $B^k (k\neq 0)$ of $B$ commutes
with $O$.} Then we let $F = f^k, F_n = T_n F T_n^{-1}$ and
consider the sequence of maps $h_n = F_n^{-1} g_n'^{-1} F_n g_n'$.
The estimate (2) of Proposition \ref{almost affine zoom} and the
hypotheses $B^k$ commutes with $O$, $B^k \neq I$ then imply that
$h_n = id + A_n + E_n$ where $A_n \in $ {\it aff}$(\mathbb{R}^N)$
tends to $0$ and $E_n / ||A_n|| \to 0$ uniformly on compacts, so
by Proposition \ref{Euler's formula} and the previous Lemma
\ref{cocompact action} we get the required flow in
$\overline{<G,f>}$.

\medskip

\noindent {\it Case 2: No non-trivial power of $B$ commutes with
$O$.} We let $f_n = T_n f T_n^{-1}$. Since $t_n' \ll
||\epsilon_n'||/\lambda_n \ll 1$, it follows that $f_n \to Df(0) =
B$ uniformly on compacts. Also by the estimate (1) of Proposition
\ref{almost affine zoom}, $g_n'^{-1}f_n g_n' \to O^{-1}BO = C$
say. By Lemma \ref{cocompact action} it suffices to show that the
closed subgroup of $GL_{N-1}(\mathbb{R})$ generated by $B,C$
contains a flow, or equivalently, that it is not discrete (since a
closed subgroup of a Lie group is a Lie group). By hypothesis $B$
is of the form $B = A tO' A^{-1}$ where $t$ is a scalar and $O'$
is orthogonal. Since no power of $B$ commutes with $O$, the linear
maps $B^{-k}C^k= B^{-k} O^{-1} B^k O = A O'^{-k}A^{-1}O^{-1}A O'^k
A^{-1} O , k \in \mathbb{Z}$ are all distinct, and moreover they
have an accumulation point since the powers of $O'$ converge along
a subsequence. $\diamond$

\medskip

We have also the following:

\medskip

\begin{prop} \label{based linear maps flow}{ Let $G$ be a cocompact Kleinian
group and $F_n$ a sequence of linear maps based at points
$\epsilon_n,M_n$. Suppose that the sequences $\epsilon_n,M_n$
contain infinitely many distinct elements, are convergent, and the
multipliers of the $F_n$'s converge to an invertible linear map of
norm less than $1$. Then the group $\hat{G} = \overline{<G,
\{F_n\}>}$ contains a non-trivial one-parameter subgroup.}
\end{prop}

\medskip

\noindent {\bf Proof:} Conjugating $G$ and the maps $F_n$ by a conformal
map sending the limits of $\epsilon_n,M_n$ to $0,\infty$ we may assume
that $\epsilon_n \to 0, M_n \to \infty$. Each $F_n$ is of the form
$S_n^{-1} A_n S_n$ where $A_n$ is a linear map and $S_n$ is the unique
conformal map such that $S_n(\epsilon_n) = 0, S_n(M_n) = \infty,
DS_n(a_n) = I$. By hypothesis, $A_n$ converges to an invertible
linear map $A$ and $S_n, S_n^{-1}$ converge to the identity,
therefore $A \in \hat{G}$.

\medskip

Now we want to adjust the poles $\epsilon_n, M_n$ of the maps
$F_n$ to get almost-affine maps as in Section 2.3. First we may
assume by passing to a subsequence if necessary that $\epsilon_n
\neq 0, M_n \neq \infty$ for all $n$. Since $\epsilon_n \to 0, M_n
\to \infty$, we can choose $t_n > 0$ such that $\epsilon_n' =
\epsilon_n / t_n, M_n' = M_n / t_n$ satisfy $1 / ||M_n'|| \ll
\epsilon_n' \ll 1$. Then putting $T_n(x) = x/t_n$, computations
very similar to those in the proof of Proposition \ref{almost
affine maps} give $F_n'(x) := T_n F_n T_n^{-1}(x) = A_n x + (I -
A_n)\epsilon_n' + o(\epsilon_n')$ uniformly on compacts. Also $A =
T_n A T_n^{-1} \in T_n \hat{G} T_n^{-1}$, so the maps $f_n =
A^{-1} F_n'$ belong to $T_n \hat{G} T_n^{-1}$, and clearly satisfy
the hypotheses of Proposition \ref{Euler's formula}. Applying
Proposition \ref{Euler's formula} and Lemma \ref{cocompact action}
gives the required one-parameter subgroup in $\hat{G}$. $\diamond$

\section{Applications}

\subsection{Pole-preserving and pole-pairing maps}

For $N \geq 3$ let $G_1, G_2$ be cocompact Kleinian groups acting
on $\partial \mathbb{H}^N$ and $h$ a $C^2$ diffeomorphism of
$\partial \mathbb{H}^N$. We consider the closed subgroup $\hat{G}
:= \overline{<G_1, h^{-1}G_2 h>}$ of Homeo$(\partial
\mathbb{H}^N)$ generated by $G_1$ and $h^{-1}G_2 h$. We denote by
Poles$(G_i), i=1,2$ the set of fixed points (or poles) of elements
of $G_i,i=1,2$. As an immediate corollary of Theorem
\ref{conformalisable}, we have

\medskip

\begin{prop} \label{C2 pole preserving}{ If the group $\hat{G}$ does not
contain a non-trivial one parameter subgroup (in particular if it
is discrete) then $h$ 'preserves poles', i.e. $h($Poles$(G_1)$) =
Poles$(G_2)$.}
\end{prop}

\medskip

\noindent {\bf Proof:} For any pole $y_0$ of an element $g_2 \in
G_2$, $x_0 = h^{-1}(y_0)$ is a fixed point of the conjugate $f =
h^{-1}g_2 h$, which by Theorem \ref{conformalisable} must be a
pole of $G_1$ under the hypothesis on $\hat{G}$. $\diamond$

\medskip

If we assume that the conformal distortion of $h$ is constant then
we can prove the stronger statement:

\medskip

\begin{theorem} \label{linear pole-pairing}{ If $h$ is equal to a linear map pre
and post-composed with conformal maps, and $\hat{G}$ does not
contain a non-trivial one parameter subgroup, then $h$ 'pairs
poles', i.e. for any $g_1 \in G_1$ with poles $p,q$ there exists
$g_2 \in G_2$ with poles $h(p), h(q)$.}
\end{theorem}

\medskip

\noindent {\bf Proof:} If $h = \phi A \psi$ with $\phi, \psi$
conformal and $A$ linear, then $<G_1, h^{-1}G_2 h> =
\psi^{-1} < (\psi G_1 \psi^{-1}), A^{-1} (\phi^{-1} G_2 \phi) A >
\psi$ so replacing $G_1, G_2$ by appropriate conjugates we may
assume that $h = A$ is linear to start with.

\medskip

By the previous Proposition $A$ preserves poles; if $h$ does not
pair poles, there are elements $g_1 \in G_1, g_2 \in G_2$ such
that $g_1$ and $g_2' = h^{-1} g_2 h$ have exactly one fixed point
$p$ in common. We proceed to make some computations to show that
some conformal conjugate of $\hat{G}$ contains a map satisfying
the hypotheses of Theorem \ref{tangent identity}.

\medskip

Conjugating by translations and replacing $h = A$
by an affine linear map $h(x) = Ax + c$ if
necessary, we may assume $p = 0, h(p) = 0$. Since $g_2$ has
a pole at $h(p) = 0$, the other pole is nonzero, so we can write
$g_2$ as the conjugate by the inversion in the unit sphere around the
origin $i(x) := x / ||x||^2$ of an affine linear map $x \mapsto \lambda O (x - \epsilon)
+ \epsilon$ with $\epsilon \neq 0, \lambda \neq 1$ and $O$ orthogonal.
Then in a neighbourhood of $x = h(p) = 0$,
\begin{align*}
g_2(x) & = \frac{\frac{\lambda O x}{||x||^2||} + (I - \lambda O)\epsilon}
{||\frac{\lambda O x}{||x||^2} + (I - \lambda O)\epsilon||^2} \\
       & = (\lambda O x + (I - \lambda O)\epsilon
       ||x||^2)\lambda^{-2}\left(1 + \frac{2<\lambda O x, (I - \lambda O)\epsilon>}{\lambda^2} +
       O(||x||^2)\right)^{-1} \\
       & = \frac{1}{\lambda}Ox + \frac{2<Ox, (\lambda O -
       I)\epsilon>}{\lambda^2}Ox + \frac{(I - \lambda O)\epsilon
       ||x||^2}{\lambda^2} + O(||x||^3) \\
\end{align*}
Since $h$ has constant derivative equal to $A$,
letting $\tilde{O} = A^{-1}OA$ and $b = (I - \lambda O)\epsilon$,
this gives, for $x$ near $p = 0$,
$$
g_2'(x) = h^{-1} g_2 h (x) = \frac{1}{\lambda}\tilde{O}x +
\frac{1}{\lambda^2}(2<A\tilde{O}x, -b>\tilde{O}x + A^{-1}b||Ax||^2) +
O(||x||^3)
$$
Now we consider the conjugate $g_2'' := i g_2' i$ of $g_2'$ by the
inversion $i$. Straightforward computations lead
to the following in a neighbourhood of $w = \infty$:
$$
g_2''(w) = \left(\lambda \tilde{O}w + \psi(w) +
O\left(\frac{1}{||w||}\right)\right)\frac{||w||^2}{||\tilde{O}w||^2}
$$
where
$$
\psi(w) = \frac{2<A\tilde{O}w, b>\tilde{O}w}{||w||^2} +
\left(A^{-1}b - \frac{2<\tilde{O}w,
A^{-1}b>\tilde{O}w}{||\tilde{O}w||^2}\right)\frac{||Aw||^2}{||w||^2}
$$
The inversion $i$ conjugates $g_1$ to an affine linear map $g_1'$ with
poles $a, \infty$ say. Now we conjugate by the translation $T(w) =
w - a$ to get a conformal linear map $g_1'' = T^{-1}g_1'T$ and a map
$g_2''' = T^{-1}g_2''T$. We note that $\psi(w + a) = \psi(w) +
O(1/||w||)$. Also, using bilinearity of the inner product and the geometric
series,
$$
\frac{||w + a||^2}{||\tilde{O}(w + a)||^2} =
\frac{||w||^2}{||\tilde{O}w||^2}\left(1 +
2\left(\frac{<w,a>}{||w||^2} - \frac{<\tilde{O}w,
\tilde{O}a>}{||\tilde{O}w||^2}\right) +
O\left(\frac{1}{||w||^2}\right)\right)
$$
Let
$$
\sigma(v, w) = \frac{<v,w>}{||w||^2} - \frac{<\tilde{O}v,
\tilde{O}w>}{||\tilde{O}w||^2}
$$
Then some calculations give
\begin{align*}
g_2'''(w) & = \left(\lambda \tilde{O}(w + a) + \psi(w + a) +
O\left(\frac{1}{||w||}\right)\right)\frac{||w + a||^2}{||\tilde{O}(w + a)||^2} -
a \\
          & = \left(\lambda \tilde{O}w + \phi(w) +
          \psi(w)\right)\frac{||w||^2}{||\tilde{O}w||^2} + O\left(\frac{1}{||w||}\right)\\
\end{align*}
where $$ \phi(w) = 2 \sigma(w,a) \lambda \tilde{O} w +
\left(\lambda \tilde{O} - \frac{||w||^2}{||\tilde{O}w||^2}\right)a
$$
Now applying Proposition \ref{zoom-in at pole} to the map $g_2'''$ and the conformal
linear map $g_1''$ we get $\hat{g_2}$, a linear map based at
$\infty, 0$ with multiplier $(1/\lambda)\tilde{O}$. The inverse
has the form $$ \hat{g_2}^{-1}(w) =
\frac{1}{\lambda}\tilde{O}^{-1}w
\frac{||w||^2}{||\tilde{O}^{-1}w||^2} $$ We consider the map $f(w)
= \hat{g_2}^{-1}(g_2'''(w))$. Let $\eta(w) = (1/\lambda)
\tilde{O}^{-1}(\phi(w) + \psi(w))$. We compute to find that
\begin{align*}
f(w) & = \left(w + \eta(w) + O\left(\frac{1}{||w||}\right)\right)(1 -
2\sigma(\eta(w),w)) + O\left(\frac{1}{||w||}\right) \\
     & = w + \Phi(w) + O\left(\frac{1}{||w||}\right) \\
\end{align*}
where
$$
\Phi(w) = \eta(w) - 2\sigma(\eta(w),w)w
$$
Since $\eta(w)$ is homogeneous of degree $0$ in $w$ and $\sigma(v,w)$
of degrees $-1$ in $v,w$, $\Phi$ is homogeneous of degree $0$. It's
straightforward to see that $\Phi(w + v) = \Phi(w) +
O(||v||/||w||)$. In order to apply Theorem \ref{tangent identity}, it remains to
check that $\Phi$ is not identically zero.

\medskip

Note that $(||\epsilon||^2/||A^{-1}\epsilon||^2)A^{-1}\epsilon$
and $a$ are the fixed points of $g_2''$ and $g_1'$ respectively which are not equal.
Since $b = (I - \lambda O)\epsilon$, replacing $g_2$ by a sufficiently large power of
$g_2$ if necessary, we may assume that $\lambda < 1$ is small
enough so that $(||\epsilon||^2/||A^{-1}b||^2)A^{-1}b \neq a$.
Now
\begin{align*}
\eta(w) = & 2w\left(\frac{<A\tilde{O}w, b>}{\lambda ||w||^2} -
\frac{||Aw||^2}{||w||^2}\frac{<\tilde{O}w, A^{-1}b>}{\lambda ||\tilde{O}w||^2}
+ \sigma(a,w)\right) \\ & + \frac{||Aw||^2}{||w||^2}
\frac{\tilde{O}^{-1}A^{-1}b}{\lambda} + \left(I -
\frac{||\tilde{O}w||^2\tilde{O}^{-1}}{\lambda||w||^2}\right)a \\
\end{align*}
We may assume $\lambda$ is small enough so that the terms involving $1/\lambda$ above are
much larger than the others. If $b = 0$, then $a \neq 0$ and $\eta(a) = t a$ for some $t \neq 0$,
so $\Phi(a) = ta - \sigma(ta, a)a = ta \neq 0$ (since $\sigma(ta,a) = 0$) and we are done.

\medskip

Otherwise if $b \neq 0$, then for $w = w_0 = \tilde{O}^{-1} A^{-1}b$,
$\eta(w_0)$ is approximately
$$
\frac{||A^{-1}b||^2}{\lambda ||\tilde{O}^{-1}
A^{-1}b||^2}\tilde{O}^{-1}\left(\frac{||b||^2}{||A^{-1}b||^2}A^{-1}b
- a\right)
$$
So if $A^{-1}b$ and $a$ are linearly independent then, assuming
$\lambda$ is small enough, $\eta(w_0)$ and
$w_0$ are linearly independent and so $\Phi(w_0) \neq 0$. Finally,
if $A^{-1}b$ and $a$ are linearly dependent, then $\eta(w_0) \neq 0$ and
$w_0$ are almost parallel in which case $\sigma(\eta(w_0),
w_0)w_0$ is small compared to $\eta(w_0)$, therefore $\Phi(w_0) \neq
0$. $\diamond$

\medskip

\noindent {\bf Proof of Theorem \ref{pole-pairing}:} It only remains to show that
the length spectra of $G_1, G_2$ are commensurable. Given $g_2 \in
G_2$, by Proposition 4.1, the map $f = h^{-1} g_2 h$ has a fixed point
$\epsilon$ in common with some $g_1 \in G_1$. By Proposition \ref{zoom-in at pole},
we get $f^* \in \hat{G}$ a linear map based at the poles of $g_1$.
So $\hat{G}$ contains the conjugates by a conformal map of two
linear maps $S = \lambda_1 O_1$ and $T = \lambda_2 O_2'$, where
$\lambda_1 O_1 = Dg_1(\epsilon), \lambda_2 O_2' = Df(\epsilon)$,
where $O_2'$ is a linear conjugate of an orthogonal map. Considering
inverses if necessary we may assume $\lambda_i < 1, \lambda_i = e^{-l_i},
i = 1,2$ where $l_i$ is the translation length of $g_i$. Since
a closed subgroup of a Lie group is a Lie group, the group $\overline{<S,T>}$
must be discrete. For $n \geq 1$, choose $m_n$ such that
$\lambda_2 < \lambda_1^n/\lambda_2^{m_n} \leq 1$. Then the sequence
of maps $S^n T^{-m_n} = (\lambda_1^n/\lambda_2^{m_n})O_1^n O_2'^{-m_n}$
contains a convergent subsequence. By discreteness this
subsequence must be eventually constant, so for some integers
$j,k$,
$$
\lambda_1^j \lambda_2^k O_1^j = O_2'^{-k}
$$
Since $O_2'$ is conjugate to an orthogonal map, the norms of
powers of the RHS above are uniformly bounded, so we must have
$\lambda_1^j \lambda_2^k = 1$, and therefore $l_1, l_2$ are rational
multiples of each other. $\diamond$

\bigskip

\subsection{Cocompact Kleinian groups and quasiconformal
maps}

Let $G_1, G_2$ be cocompact Kleinian groups in dimension $N \geq
3$, $h$ a quasiconformal map of $\partial \mathbb{H}^N$ and
$\hat{G} = \overline{<G_1, h^{-1} G_2 h>}$.

\medskip

\begin{lemma} \label{upgrade to linear}{ There are cocompact Kleinian
groups $G^*_1, G^*_2$ and a linear map $A$ such that $\hat{G}$
contains a conformal conjugate
of $\overline{<G^*_1, A^{-1} G^*_2 A>}$.
If $h$ is not conformal, $A$ can be taken to be non-conformal.}
\end{lemma}

\medskip

\noindent {\bf Proof:} We recall that quasiconformal maps are
differentiable almost everywhere. Let $x_0$ be a point of
differentiability of $h$. We note if $h$ is not conformal
then we can choose $x_0$ such that $Dh(x_0)$ is not conformal.
Given conformal maps $\phi_1,\phi_2$,
we note $\hat{G} = \phi_1 \overline{< \phi_1^{-1} G_1 \phi_1,
(\phi_2 h \phi_1)^{-1} (\phi_2 G_2 \phi_2^{-1}) (\phi_2 h
\phi_1)>} \phi_1^{-1}$, so replacing $G_1, G_2$ by conformal
conjugates if necessary we may compose $h$ with conformal maps
on the left and right to assume that $x_0 = h(x_0) = 0$.

\medskip

Let $T_n(x) = x/n$. The maps $T_n^{-1} h T_n$ converge to the
linear map $Dh(0) = A$ say. By Lemma \ref{cocompact action} we may write $T_n =
g_{1,n} \phi_n = g_{2,n} \psi_n$ where $g_{i,n} \in G_i, i=1,2$
and $\phi_n, \psi_n$ are conformal maps converging along
subsequences to conformal maps $\phi, \psi$ respectively. Then the
maps $g_{2,n}^{-1} h g_{1,n} = \psi_n (T_n^{-1} h T_n)
\phi_n^{-1}$ converge to $\psi A \phi^{-1}$.

\medskip

Thus for any $g_2 \in G_2$, the maps
$$
(g_{2,n}^{-1} h g_{1,n})^{-1} g_2 (g_{2,n}^{-1} h g_{1,n}) =
g_{1,n}^{-1}(h^{-1}(g_{2,n}g_2 g_{2,n}^{-1})h)g_{1,n} \in \hat{G}
$$
converge to a map $\phi A^{-1} \psi^{-1} g_2 \psi A \phi^{-1}$.
It follows that $\hat{G}$ contains the group $\phi (A^{-1} G_2^*
A) \phi^{-1}$ where $G_2^* = \psi^{-1} G_2 \psi$. Letting $G_1^* = \phi^{-1}
G_1 \phi$, since $\hat{G}$ also contains $\phi G_1^* \phi^{-1} =
G_1$, it follows that $\hat{G}$ contains $\phi \overline{<G_1^*,
A^{-1}G_2^* A>} \phi^{-1}$. $\diamond$

\medskip

\noindent {\bf Proof of Theorem \ref{main theorem}:} By the previous Lemma, we
may assume that $h = A$ is a non-conformal linear map. Suppose
$\hat{G}$ does not contain a non-trivial one-parameter subgroup.
Then by Theorem \ref{linear pole-pairing}, the linear map $A$ pairs poles of $G_1$ with
those of $G_2$.

\medskip

Pick $g_1 \in G_1$ with both poles distinct from $0,\infty$, and suppose $A$ pairs
these poles with those of some $g_2 \in G_2$. Let $\phi_1, \phi_2$ be
conformal maps sending the poles of $g_1,g_2$ respectively to $0,
\infty$. The maps $g_j' := \phi_j g_j \phi_j^{-1} \in G_j'
:= \phi_j G_j \phi_j^{-1}$ are
conformal linear maps. The smooth map $\mu := \phi_2 A \phi_1^{-1}$
pairs poles of $G_1', G_2'$ (because $A$ pairs poles of $G_1,G_2$),
satisfies $\mu(0) = 0, \mu(\infty) = \infty$, and the group
$\hat{G}' := \phi_1^{-1} \overline{< G_1', \mu^{-1} G_2' \mu>} \phi_1 =
\phi_1 \hat{G} \phi_1^{-1}$ does not contain a non-trivial
one-parameter subgroup.
Moreover, because $A$ is a non-conformal linear map and $\phi_j$'s
are non-linear conformal maps, it is not hard to see that $\mu$ is
not a linear map.

\medskip

Let $g_1' = \lambda_1 O_1, g_2' = \lambda_2 O_2$ where the $O_j$'s are
orthogonal, and $\lambda_j$'s positive scalars; we may assume
$\lambda_j < 1$. For every $n \geq 1$ let $m_n \geq 1$ be such that
$\lambda_2 \leq \lambda_1^n / \lambda_2^{m_n} \leq 1$. Since $\mu(0) = 0$
and $\mu$ is differentiable at $0$, it is
clear that along a subsequence the maps
$$
h_n(x) := g_2'^{-m_n} \mu g_1'^n (x) = \left(\frac{\lambda_1^n}{\lambda_2^{m_n}}\right)O_2^{-m_n}
\left(\frac{1}{\lambda_1^n}\mu(\lambda_1^n O_1^n x)\right)
$$
and their inverses converge uniformly on compacts to linear maps $B, B^{-1}$ say.
It is not hard to see that in fact, since $\mu$ is smooth,
we have $C^1$ convergence on compacts. By an argument similar
to the proof of Lemma \ref{upgrade to linear}, the group $\overline{< G_1',
B^{-1} G_2' B>}$ is contained in $\hat{G}'$, hence by Theorem \ref{linear pole-pairing}, the linear
map $B$ must pair poles of $G_1'$ and $G_2'$.
Since none of the maps $h_n$ are linear, there are infinitely
many distinct maps in the sequence $h_n$.

\medskip

In the terminology of Schwartz \cite{schwarz-inv}, the maps $h_n$ are {\it eccentric
maps}. The following lemma is due to Schwartz:

\medskip

\begin{lemma} \label{scattering} (\cite{schwarz-inv}){ For any open set $U \subset
\mathbb{R}^N$ there exists a constant $\delta = \delta(U)$ such
that if $\Sigma \subset U$ is $\delta$-dense in $U$ (i.e. $U$ is
contained in the $\delta$-neighbourhood of $\Sigma$) then any two
eccentric maps which agree on $\Sigma$ agree everywhere.}
\end{lemma}

\medskip

Since poles of $G_2'$ are dense and there are infinitely many distinct maps
$h_n$, it follows from the above lemma that there is a $g \in G_2'$ with poles
$a,b$ such that $a_n := h_n^{-1}(a) \neq B^{-1}a, b_n := h_n^{-1}(b) \neq B^{-1}b$ for
infinitely many $n$. Since $h_n$ pairs poles of $G_1', G_2'$, for all $n$
there is a $g_n \in G_1'$ with
poles $a_n, b_n$. Applying Proposition \ref{zoom-in at pole} to the maps $f_n =
h_n^{-1} g h_n \in \hat{G}'$ and $g_n$, we get $F_n \in \hat{G}'$, a linear map based at
$a_n,b_n$. The $C^1$ convergence of the $h_n$'s
implies that the multipliers of the $F_n$'s converge. Applying
Proposition \ref{based linear maps flow} to the sequence $F_n$ gives a non-trivial
one-parameter group of maps in $\hat{G}'$,
contradicting our initial hypothesis. $\diamond$

\bigskip

\subsection{Pattern rigidity}

\bigskip

Let $\mathcal{J}_1,
\mathcal{J}_2$ be collections of closed subsets of $\partial \mathbb{H}^N$
invariant under cocompact Kleinian groups $G_1, G_2$
respectively and discrete (with respect to the Hausdorff topology)
in the space of compact subsets minus singletons.

\medskip

\begin{lemma} \label{discrete no flow} { A collection of closed
subsets $\mathcal{J}$ of $\partial \mathbb{H}^N$ which is discrete
in the space of compact subsets minus singletons cannot be invariant under a non-trivial
continuous one-parameter group
$(f_t)_{t \in \mathbb{R}}$ of homeomorphisms of $\partial \mathbb{H}^N$.}
\end{lemma}

\medskip

\noindent {\bf Proof:} If there is such a one-parameter group,
let $t_0 \in \mathbb{R}, x \in \partial
\mathbb{H}^N$ be such that $f_{t_0}(x) \neq x$. We can find a
$J \in \mathcal{J}$ contained in a sufficiently small
neighbourhood of $x$ such that $f_{t_0}(J) \cap J = \emptyset$.
On the other hand since $f_t \to id$ as $t \to 0$, by discreteness
of $\mathcal{J}$ with respect to the Hausdorff topology, we must
have $f_{t_0/n}(J) = J$ for $n$ large enough, so $f_{t_0}(J)
= f_{t_0/n}^n(J) = J$, a contradiction. $\diamond$

\medskip

Pattern rigidity is now an almost immediate corollary of Theorem
\ref{main theorem}:

\medskip

\noindent{\bf Proof of Corollary \ref{pattern rigidity}:} Discreteness in the Hausdorff
topology implies that any uniform limit of maps pairing $\mathcal{J}_1$ with $\mathcal{J}_2$
must also pair $\mathcal{J}_1$ with $\mathcal{J}_2$.
In particular the group of homeomorphisms preserving a discrete collection
$\mathcal{J}$ is a closed subgroup of Homeo$(\partial \mathbb{H}^N)$.
Let $h$ be a quasiconformal
map pairing $\mathcal{J}_1$ with $\mathcal{J}_2$. Then the
groups $G_1, h^{-1} G_2 h$ preserve the collection $\mathcal{J}_1$
and hence so does the group $\hat{G} = \overline{<G_1, h^{-1} G_2 h>}$.
By Lemma \ref{discrete no flow}, $\hat{G}$ has no non-trivial one-parameter subgroups,
so it follows from Theorem \ref{main theorem} that $h$ must be conformal.$\diamond$

\bibliography{qcflows}
\bibliographystyle{alpha}

\medskip

Ramakrishna Mission Vivekananda University,
Belur Math, WB-711202, India

\end{document}